\documentclass[10pt]{article}
\usepackage{amssymb, amsmath, url, graphicx, setspace, geometry}
\usepackage{calrsfs}
\usepackage{wasysym}

\def\3{\subset }
\def\4{\subseteq }
\def\<{\left<}
\def\>{\right>}

\def\bit{\begin{itemize}}
\def\eit{\end{itemize}}
\def\3{\subset }
\def\4{\subseteq }

\def\0{\leqno}

\def\barr{\begin{array}}
\def\earr{\end{array}}

\def\Z{{\rlap{$\kern2pt{\rm Z}$}{\rm Z}\,}}
\def\bld#1#2{{\buildrel{#1}\over{#2}}}
\def\st#1#2{{\mathrel{\mathop{#2}\limits_{#1}}{}\!}}
\def\stb#1#2#3{{\st{{#1}}{\bld{{#2}}{#3}}{}\!}}
\def\xmare#1#2{\stb{#1}{#2}{\mbox{\Huge$\times$}}}

\title{\bf On a divisibility property involving the sum of element orders}
\author{Mihai-Silviu Lazorec\footnote{The author was supported by the European Social Fund, through Operational Programme Human Capital 2014-2020, project no. POCU/380/6/13/123623.}}
\date{March 3, 2020}

\begin{document}
\maketitle

\begin{abstract}
A finite group $G$ is called $\psi$-divisible if $\psi(H)|\psi(G)$ for any subgroup $H$ of $G$, where $\psi(H)$ and $\psi(G)$ are the sum of element orders of $H$ and $G$, respectively. In this paper, we extend a result provided in \cite{10}, by classifying the finite groups whose all subgroups are $\psi$-divisible. Since the existence of $\psi$-divisible groups is related to the class of square-free order groups, we also study the sum of element orders and the $\psi$-divisibility property of ZM-groups. In the end, we introduce the concept of $\psi$-normal divisible group, i.e. a group for which the $\psi$-divisibility property is satisfied by all its normal subgroups. Using simple and quasisimple groups, we are able to construct infinitely many $\psi$-normal divisible groups which are neither simple nor nilpotent.   
\end{abstract}



\section{Introduction}

Let $G$ be a finite group. For an element $x$ of $G$, we denote the order of $x$ by $o(x)$. In \cite{1}, Amiri, Jafarian Amiri and Isaacs, introduced the sum of element orders of $G$, which is a quantity denoted by $\psi(G)$ and, as its name suggests, it is defined as 
$$\psi(G)=\sum\limits_{x\in G}o(x).$$
Their main result states that the maximum value of the sum of element orders among all groups of order $|G|$ is attained if one works with the cyclic group $\mathbb{Z}_{|G|}$, i.e.  $\psi(G)\leq \psi(\mathbb{Z}_{|G|})$ and the equality holds if and only if $G\cong \mathbb{Z}_{|G|}$. After \cite{1} was published, a lot of interesting papers used the sum of element orders to characterize the nature (abelian, nilpotent, supersolvable, solvable) of a finite group $G$ via some inequalities of the following type: $\psi(G)> \alpha\psi(\mathbb{Z}_{|G|}),$, where $\alpha\in (0,1)$ is a constant. For more details, we refer the reader to \cite{4, 5, 11, 12, 15}. 

Instead of using such inequalities, in \cite{10}, Harrington, Jones and Lamarche choose to work with the usual divisibility property of positive integers. They say that a finite group $G$ is \textit{$\psi$-divisible} if 
$$\psi(H)|\psi(G), \ \forall \ H\in L(G),$$
where $G$ is the subgroup lattice of $G$. In what follows, we recall their main result which refers to a characterization of cyclic groups of square-free order via the $\psi$-divisibility property. \\

\textbf{Theorem 1.1.} \textit{Let $G$ be a finite abelian group. Then $G$ is $\psi$-divisible if and only if $G$ is cyclic of square-free order.}\\

In other words, the only $\psi$-divisible abelian groups are the cyclic ones of square-free order. 
In the same paper (see the paragraph ``3.1. A Remark on Nonabelian $\psi$-Divisible Groups"), the authors say that they were not able to find a non-abelian $\psi$-divisible group and they show that there are infinitely non-abelian groups which are not $\psi$-divisible. As a consequence of Theorem 1.1 (see Corollary 16 of \cite{10}), it is clear that given a finite abelian group $G$, all its subgroups are $\psi$-divisible if and only if $G$ is cyclic of square-free order. In the second section of this paper,  
we generalize this last statement by proving that all subgroups of a finite group $G$ are $\psi$-divisible if and only if $G$ is cyclic of square-free order. This generalization and Theorem 1.1 outline the fact that the square-free order and the $\psi$-divisibility properties are quite related. Since all square-free order groups are Zassenhaus metacyclic groups (ZM-groups, in short), we highlight some results on the sum of element orders of ZM-groups. They constitute additional arguments for the possible nonexistence of non-abelian $\psi$-divisible groups. We mention that, in \cite{13}, Jafarian Amiri and Amiri also studied some connections between the sum of element orders and the groups of square-free order.

In the third section, instead of working with all subgroups of a finite group $G$, we restrict our study by choosing to focus only on its normal subgroups. We say that $G$ is a \textit{$\psi$-normal divisible} group if $$\psi(H)|\psi(G), \ \forall \ H\in N(G),$$
where $N(G)$ is the normal subgroup lattice of $G$. Besides the cyclic groups of square-free order, we are able to indicate a lot more examples of $\psi$-normal divisible groups, most of them being constructed based on simple and quasisimple groups. 
 
Some open problems are highlighted throughout the paper. Most of our notation is standard and will usually not be repeated here. Elementary notions and results on groups can be found in \cite{9}.

\section{On the existence of non-nilpotent $\psi$-divisible groups}

We begin our study by recalling some preliminary results on the sum of element orders of a finite group.\\
 
\textbf{Lemma 2.1.} \textit{The following statements hold:
\begin{itemize}
\item[i)] (see Lemma 2.9 (1) of \cite{11}) $\psi(\mathbb{Z}_{p^n})=\frac{p^{2n+1}+1}{p+1}$, where $p$ is a prime and $n\geq 1$ is an integer. 
\item[ii)] (see Lemma 2.1 of \cite{2}) $\psi$ is multiplicative, i.e. if $G_1$ and $G_2$ are two finite groups such that $(|G_1|, |G_2|)=1$ then $\psi(G_1\times G_2)=\psi(G_1)\psi(G_2)$.
\item[iii)] (see Lemma 2.2 (5) of \cite{11}) Let $G=P\rtimes H$, where $P$ is a cyclic $p$-group, $|H|>1$ and $(p, |H|)=1$. Then,
$$\psi(G)=|P|\psi(H)+(\psi(P)-|P|)\psi(C_H(P)).$$
\end{itemize}}


Also, we recall some ideas on ZM-groups by following Theorem 9.4.3 of \cite{9}. These groups are exactly the ones for which all Sylow subgroups are cyclic. Their structure is given as follows
$$ZM(m,n,r)=\langle a,b \ | \ a^m=b^n=1, \ b^{-1}ab=a^r\rangle,$$
where the triple $(m,n,r)\in\mathbb{N}^3$ is chosen such that $(m,n)=(m,r-1)=1$ and $r^n\equiv 1 \ (mod \ m)$. The order of $ZM(m,n,r)$ is $mn$ and, using the previous arithmetic properties of $m, n$ and $r$, one can easily show that $m$ is odd. To make computations in $ZM(m,n,r)$, it is helpful to know that:
$$\alpha(x_1,y_1)\alpha(x_2,y_2)=\alpha(x_1+x_2,r^{x_2}y_1+y_2),$$
where $\alpha(x,y)=b^xa^y$, for all $(x,y)\in \lbrace 0,1,\ldots, n-1\rbrace\times\lbrace 0,1,\ldots, m-1\rbrace$. Finally, we mention that the subgroups of a ZM-group were determined in \cite{6} (see the proof of the main result). Note that if $m\ne 1$ and $n\ne 1$, the corresponding group $ZM(m,n,r)$ is non-nilpotent since it contains at least one non-normal maximal subgroup.   

In what follows, we are interested in determining a necessary and sufficient condition such that the $\psi$-divisibility property of a finite group is inherited by its subgroups.\\   

\textbf{Theorem 2.2.} \textit{Let $G$ be a finite group. Then all subgroups of $G$ are $\psi$-divisible if and only if $G$ is cyclic of square-free order.}\\

\textbf{Proof.} Let $G$ be a finite group such that all its subgroups are $\psi$-divisible. Assume that $G$ is not of square-free order. Then, there is a prime number $p$ such that $p^2||G|$, so $G$ has a subgroup $H$ of order $p^2$. Hence, $H$ is abelian and, according to our hypothesis, it is $\psi$-divisible, contradicting the classification provided by Theorem 1.1.  

Since $G$ is of square-free order, all its Sylow subgroups are cyclic. Then $G\cong ZM(m,n,r)$, where $m$ and $n$ are square-free numbers. We recall that $m$ must be odd. If $m=1$, then $G\cong \mathbb{Z}_n$, and we are done. If $m\geq 3$, let $p_1, p_2,\ldots, p_k$ be prime numbers such that $|G|=mn=p_1p_2\ldots p_k$. To complete our proof, we must show that $n=1$, i.e. $G\cong \mathbb{Z}_m$. Assume that $n\ne 1$. Then, there exists a pair $(i,j)\in\lbrace 1,2,\ldots, k\rbrace\times \lbrace 1,2,\ldots, k\rbrace$, with $i\ne j$, such that $p_i|n$, $p_j|m$ and $p_i|p_j-1$. We consider a subgroup $H\cong \mathbb{Z}_{p_j}\rtimes \mathbb{Z}_{p_i}$ of $G$. Since $H$ is $\psi$-divisible, we have
$$\psi(\mathbb{Z}_{p_j})|\psi(H)\Longleftrightarrow p_j^2-p_j+1|p_j^2-p_j+1+(p_i-1)p_ip_j,$$ which further implies that $$p_j^2-p_j+1|(p_i-1)p_ip_j.$$
Since $p_i, p_j$ are primes and $p_i|p_j-1$, it is easy to show that $(p_j^2-p_j+1,p_j)=(p_j^2-p_j+1,p_i)=1.$ Consequently $p_j^2-p_j+1|p_i-1$, which leads to $p_j^2-p_j+1\leq p_i-1$. But $p_i\leq p_j-1$, so $p_j^2-2p_j+3\leq 0$, a contradiction. Hence, $n=1$ and $G\cong \mathbb{Z}_m.$ 

Since the converse is true due to Theorem 1.1, our proof is complete. 
\hfill\rule{1,5mm}{1,5mm}\\

According to Theorem 2.2, it seems that there is a quite strong connection between the $\psi$-divisibility and the square-free order properties of finite groups. As we mentioned in our previous proof, a group of square-free order is a ZM-group. Hence, if one would want to find some non-abelian $\psi$-divisibile groups, a starting point is to check the ZM-groups (especially the ones of square-free orders). In what follows, we prove some results that provide a partial answer related to this aspect, as well as an explicit formula that can be used to determine the sum of element orders of any ZM-group.\\

\textbf{Proposition 2.3.} \textit{Let $ZM(m,n,r)$ be a ZM-group such that $m=p_1^{\alpha_1}p_2^{\alpha_2}\ldots p_k^{\alpha_k}$, where $p_i$ is a prime number and $\alpha_i\geq 1$ is an integer for all $i\in\lbrace 1,2,\ldots, k \rbrace$. Then
\begin{equation}\label{r1}
\psi(ZM(m,n,r))=m\psi(\mathbb{Z}_n)+\sum\limits_{J\subsetneq I}\prod\limits_{i\in J}p_i^{\alpha_i}\prod\limits_{i\in I\setminus J}(\psi(C_{{p_i}^{\alpha_i}})-p_i^{\alpha_i})\psi(\langle b^{\beta_{I\setminus J}\cdot o_{M_{IJ}}^{I\setminus (J\cup\lbrace M_{IJ}\rbrace)}}\rangle),
\end{equation}
where $I=\lbrace 1,2,\ldots, k\rbrace, m_{IJ}=\min\limits_{i\in I\setminus J}i, M_{IJ}=\max\limits_{i\in I\setminus J}i$ and $\beta_{I\setminus J}$ is the product of all multiplicative orders needed to compute the quantity
$$o_{M_{IJ}}^{I\setminus (J\cup\lbrace M_{IJ}\rbrace)}=\begin{cases} the \ multiplicative \ order \ of \ r \ modulo \  p^{\alpha_{m_{IJ}}}_{m_{IJ}} &\mbox{, if } |J|=|I|-1 \\ the \ multiplicative \ order \ of \ r^{o_{M_{IJ}}^{I\setminus (J\cup \lbrace m_{IJ}, M_{IJ}\rbrace)}} \ modulo \ p^{\alpha_{m_{IJ}}}_{m_{IJ}} &\mbox{, if } |J|<|I|-1\end{cases}.$$}

\textbf{Proof.} Before outlining our reasoning, we would like to highlight a quite suggestive example on how the quantities $\beta_{I\setminus J}$ are computed. Hence, for $I=\lbrace 1,2,3,4,5\rbrace$ and $J=\lbrace 2, 4\rbrace$, we would like to determine $o_{M_{IJ}}^{I\setminus (J\cup\lbrace M_{IJ}\rbrace)}=o_5^{\lbrace 1,3\rbrace}$, which is the multiplicative order of $r^{o_5^{\lbrace 3\rbrace}}$ modulo $p_1^{\alpha_1}$. For this purpose, we first need to compute $o_5^{\lbrace 3\rbrace}$, which is the multiplicative order of $r^{o_5^\emptyset}$ modulo $p_3^{\alpha_3}$. Consequently, we must determine $o_5^\emptyset$, which is the multiplicative order of $r$ modulo $p_5^{\alpha_5}$. Hence, $\beta_{\lbrace 1,3,5\rbrace}=o_5^\emptyset\cdot o_5^{\lbrace 3\rbrace}$.  

Let $ZM(m,n,r)$ be a ZM-group such that $m=p_1^{\alpha_1}p_2^{\alpha_2}\ldots p_k^{\alpha_k}$. Once can check that this group may be represented as an iterated semidirect product as follows
$$ZM(m,n,r)\cong \mathbb{Z}_{p_k^{\alpha_k}}\rtimes (\mathbb{Z}_{p_{k-1}^{\alpha_{k-1}}}\rtimes (\mathbb{Z}_{p_{k-2}^{\alpha_{k-2}}}\rtimes(\ldots \rtimes(\mathbb{Z}_{p_1^{\alpha_1}}\rtimes \mathbb{Z}_n)\ldots ))).$$
We proceed by induction on $k$. If $k=1$, then $m=p_1^{\alpha_1}$ and, by Lemma 2.1 iii), we have
$$\psi(ZM(m,n,r))=p_1^{\alpha_1}\psi(\mathbb{Z}_n)+(\psi(\mathbb{Z}_{p_1^{\alpha_1}})-p_1^{\alpha_1})\psi(C_{\mathbb{Z}_n}(\mathbb{Z}_{p_1^{\alpha_1}})).$$
Remark that
\begin{align*}
b^j\in C_{\mathbb{Z}_n}(\mathbb{Z}_{p_1^{\alpha_1}})&\Longleftrightarrow b^ja^i=a^ib^j, \ \forall \ i\in \lbrace 0,1,\ldots, p_1^{\alpha_1}-1\rbrace \\& \Longleftrightarrow b^ja^i=b^ja^{r^ji}, \ \forall \ i\in \lbrace 0,1,\ldots, p_1^{\alpha_1}-1\rbrace \\ & \Longleftrightarrow r^j\equiv 1 \ (mod \ p_1^{\alpha_1}) \\ & \Longleftrightarrow o_1^\emptyset|j,
\end{align*}
where $o_1^\emptyset$ is the multiplicative order of $r$ modulo $p_1^{\alpha_1}$. Hence $C_{\mathbb{Z}_n}(\mathbb{Z}_{p_1^{\alpha_1}})=\langle b^{o_1^\emptyset}\rangle$, and
$$\psi(ZM(m,n,r))=p_1^{\alpha_1}\psi(\mathbb{Z}_n)+(\psi(\mathbb{Z}_{p_1^{\alpha_1}})-p_1^{\alpha_1})\psi(\langle b^{o_1^\emptyset}\rangle),$$
which is exactly the same as formula (\ref{r1}) for $I=\lbrace 1\rbrace$.

Further, we suppose that (\ref{r1}) holds for $k-1$, i.e. for any ZM-group $ZM(m',n',r')$ such that $m'$ is divisible by exactly $k-1$ prime numbers. Let $G\cong \mathbb{Z}_{p_{k-1}^{\alpha_{k-1}}}\rtimes (\mathbb{Z}_{p_{k-2}^{\alpha_{k-2}}}\rtimes(\ldots \rtimes(\mathbb{Z}_{p_1^{\alpha_1}}\rtimes \mathbb{Z}_n)\ldots))$ and consider the sets $I=\lbrace 1,2,\ldots, k\rbrace, I'=I\setminus \lbrace k\rbrace$. By Lemma 2.1 iii), we have
$$\psi(ZM(m,n,r))=p_k^{\alpha_k}\psi(G)+(\psi(\mathbb{Z}_{p_k^{\alpha_k}})-p_k^{\alpha_k})\psi(C_G(\mathbb{Z}_{p_k^{\alpha_k}})).$$
Since $G$ and
$$C_G(\mathbb{Z}_{p_k^{\alpha_k}}))=\mathbb{Z}_{p_{k-1}^{\alpha_{k-1}}}\rtimes (\mathbb{Z}_{p_{k-2}^{\alpha_{k-2}}}\rtimes(\ldots \rtimes(\mathbb{Z}_{p_1^{\alpha_1}}\rtimes \langle b^{o_k^\emptyset}\rangle)\ldots ))$$
are both ZM-groups satisfying the inductive hypothesis, the sum of element orders of $ZM(m,n,r)$ is
\begin{align*}
&p_k^{\alpha_k}\bigg(\frac{m}{p_k^{\alpha_k}}\psi(\mathbb{Z}_n)+\sum\limits_{J\subsetneq I'}\prod\limits_{i\in J}p_i^{\alpha_i}\prod\limits_{i\in I'\setminus J}(\psi(\mathbb{Z}_{{p_i}^{\alpha_i}})-p_i^{\alpha_i})\psi(\langle b^{\beta_{I'\setminus J}\cdot o_{M_{I'J}}^{I'\setminus (J\cup\lbrace M_{I'J}\rbrace)}} \rangle)\bigg)\\ &+(\psi(\mathbb{Z}_{p_k^{\alpha_k}})-p_k^{\alpha_k})\bigg(\frac{m}{p_k^{\alpha_k}}\psi(\langle b^{o_k^\emptyset}\rangle)+\sum\limits_{J\subsetneq I'}\prod\limits_{i\in J}p_i^{\alpha_i}\prod\limits_{i\in I'\setminus J}(\psi(\mathbb{Z}_{{p_i}^{\alpha_i}})-p_i^{\alpha_i})\psi(\langle b^{\beta_{I'\setminus J}\cdot o_{M_{I'J}}^{I'\setminus (J\cup\lbrace M_{I'J}\rbrace)}}\rangle)\bigg)\\
&=m\psi(\mathbb{Z}_n)+\frac{m}{p_k^{\alpha_k}}(\psi(\mathbb{Z}_{p_k^{\alpha_k}})-p_k^{\alpha_k})\psi(\langle b^{o_k^\emptyset}\rangle)\\&+\sum\limits_{J\subsetneq I'}\bigg(\prod\limits_{i\in J\cup\lbrace k\rbrace}p_i^{\alpha_i}\prod\limits_{i\in I'\setminus J}(\psi(\mathbb{Z}_{{p_i}^{\alpha_i}})-p_i^{\alpha_i})\psi(\langle b^{\beta_{I'\setminus J}\cdot o_{M_{I'J}}^{I'\setminus (J\cup\lbrace M_{I'J}\rbrace)}} \rangle)\\& +\prod\limits_{i\in J}p_i^{\alpha_i}\prod\limits_{i\in (I'\cup\lbrace k\rbrace)\setminus J}(\psi(\mathbb{Z}_{{p_i}^{\alpha_i}})-p_i^{\alpha_i})\psi(\langle b^{\beta_{(I'\cup\lbrace k\rbrace)\setminus J}\cdot o_{M_{(I'\cup\lbrace k\rbrace)J}}^{(I'\cup\lbrace k\rbrace)\setminus (J\cup\lbrace M_{(I'\cup\lbrace k\rbrace)J}\rbrace)}} \rangle)\bigg)\\&
=m\psi(\mathbb{Z}_n)+\frac{m}{p_k^{\alpha_k}}(\psi(\mathbb{Z}_{p_k^{\alpha_k}})-p_k^{\alpha_k})\psi(\langle b^{o_k^\emptyset}\rangle)\\&+\sum\limits_{J\subsetneq I'}\bigg(\prod\limits_{i\in J\cup\lbrace k\rbrace}p_i^{\alpha_i}\prod\limits_{i\in I'\setminus J}(\psi(\mathbb{Z}_{{p_i}^{\alpha_i}})-p_i^{\alpha_i})\psi(\langle b^{\beta_{I'\setminus J}\cdot o_{M_{I'J}}^{I'\setminus (J\cup\lbrace M_{I'J}\rbrace)}} \rangle)\\&+ \prod\limits_{i\in J}p_i^{\alpha_i}\prod\limits_{i\in I\setminus J}(\psi(\mathbb{Z}_{{p_i}^{\alpha_i}})-p_i^{\alpha_i})\psi(\langle b^{\beta_{I\setminus J}\cdot o_{M_{IJ}}^{I\setminus (J\cup\lbrace M_{IJ}\rbrace)}} \rangle)\bigg)\\&
=m\psi(\mathbb{Z}_n)+\sum\limits_{J\subsetneq I}\prod\limits_{i\in J}p_i^{\alpha_i}\prod\limits_{i\in I\setminus J}(\psi(\mathbb{Z}_{{p_i}^{\alpha_i}})-p_i^{\alpha_i})\psi(\langle b^{\beta_{I\setminus J}\cdot o_{M_{IJ}}^{I\setminus (J\cup\lbrace M_{IJ}\rbrace)}} \rangle),
\end{align*}
which completes our proof.
\hfill\rule{1,5mm}{1,5mm}\\

As an application, we use (\ref{r1}) to compute the sum of element orders of $G\cong ZM(105,4,62)$ (SmallGroup(420, 23)). We have
\begin{align*}
&\psi(G)=105\psi(\mathbb{Z}_{4})+15(\psi(\mathbb{Z}_7)-7)\psi(\langle b^{o_3^\emptyset}\rangle)+21(\psi(\mathbb{Z}_5)-5)\psi(\langle b^{o_2^\emptyset}\rangle)+35(\psi(\mathbb{Z}_3)-3)\psi(\langle b^{o_1^\emptyset}\rangle)\\&+3(\psi(\mathbb{Z}_5)-5)(\psi(\mathbb{Z}_7)-7)\psi(\langle b^{o_3^\emptyset o_3^{\lbrace 2 \rbrace}}\rangle)+5(\psi(\mathbb{Z}_3)-3)(\psi(\mathbb{Z}_7)-7)\psi(\langle b^{o_3^\emptyset o_3^{\lbrace 1 \rbrace}}\rangle)\\&+7(\psi(\mathbb{Z}_3)-3)(\psi(\mathbb{Z}_5)-5)\psi(\langle b^{o_2^\emptyset o_2^{\lbrace 1 \rbrace}}\rangle)+(\psi(\mathbb{Z}_3)-3)(\psi(\mathbb{Z}_5)-5)(\psi(\mathbb{Z}_7)-7)\psi(\langle b^{o_3^\emptyset o_3^{\lbrace 2\rbrace} o_3^{\lbrace 1,2\rbrace}}\rangle).
\end{align*}
The multiplicative orders are $o_1^\emptyset=2, o_2^\emptyset=4, o_3^\emptyset=2, o_2^{\lbrace 1\rbrace}=1, o_3^{\lbrace 1\rbrace}=1, o_3^{\lbrace 2\rbrace}=2, o_3^{\lbrace 1,2\rbrace}=1$. Making the replacements in the previous formula and using Lemma 2.1 i), we obtain $\psi(G)=10171.$ The result can be checked via GAP \cite{16}.\\

Further, we show that there are infinitely many non-nilpotent ZM-groups that are not $\psi$-divisible.\\

\textbf{Proposition 2.4.} \textit{Let $ZM(m,n,r)$ be a ZM-group such that $m=p_1^{\alpha_1}$ and $n=q_1^{\beta_1}q_2^{\beta_2}\ldots q_s^{\beta_s}$, where $p_1, q_1, q_2, \ldots, q_s$ are prime numbers and $\alpha_1, \beta_1,\beta_2,\ldots, \beta_s$ are positive integers. Consider the set $A=\lbrace i\in\lbrace 1,2,\ldots, s\rbrace \ | \ q_i|p_1-1\rbrace$. If $\alpha_1\geq \max\limits_{i \in A}\beta_i$, then $ZM(m,n,r)$ is not $\psi$-divisible.}\\

\textbf{Proof.}
Let $G\cong ZM(m,n,r)$ be a ZM-group as indicated by our hypothesis. Note that $A\ne \emptyset$ since $G$ is not cyclic. According to Lemma 2.1 iii), we have
$$\psi(G)=p_1^{\alpha_1}\psi\bigg(\xmare{i=1}s \mathbb{Z}_{q_i^{\beta_i}}\bigg)+(\psi(\mathbb{Z}_{p_i^{\alpha_i}})-p_i^{\alpha_i})\psi(\langle b^{o_1^\emptyset}\rangle).$$
Since $o_1^\emptyset$ is the multiplicative order of $r$ modulo $p_1^{\alpha_1}$, it easily follows that $r^{o_1^\emptyset}\equiv 1 \ (mod \ p_1)$. Hence $p_1-1|o_1^{\emptyset}$, so all elements of $A$ appear in the prime factorization of $o_1^{\emptyset}$. Then
$$\langle b^{o_1^\emptyset}\rangle\cong G_1\times G_2, \text{ where } G_1\cong\xmare{i\not\in A}{}\mathbb{Z}_{q_i^{\beta_i}} \text{ and } G_2\cong \xmare{i\in A}{}\mathbb{Z}_{q_i^{\beta_i-\gamma_i}},$$ with $0\leq\gamma_i\leq\beta_i, \ \forall \ i\in A$, and at least one positive $\gamma_i$.

Assume that $G$ is $\psi$-divisible and consider the cyclic subgroup $H\cong\mathbb{Z}_{p_1^{\alpha_1}}\times G_1$ of $G$. Then, by applying Lemma 2.1 ii), we get
\begin{align*}
\psi(H)|\psi(G)&\Longleftrightarrow \psi(\mathbb{Z}_{p_1^{\alpha_1}})\psi(G_1)\bigg|\bigg[p_1^{\alpha_1}\psi\bigg(\xmare{i\in A}{}\mathbb{Z}_{q_i^{\beta_i}} \bigg)+(\psi(\mathbb{Z}_{p_1^{\alpha_1}})-p_1^{\alpha_1})\psi(G_2)\bigg]\psi(G_1)\\
&\Longleftrightarrow \psi(\mathbb{Z}_{p_1^{\alpha_1}})\bigg|p_1^{\alpha_1}\bigg[\psi\bigg(\xmare{i\in A}{}\mathbb{Z}_{q_i^{\beta_i}} \bigg)-\psi(G_2)\bigg]+\psi(\mathbb{Z}_{p_1^{\alpha_1}})\psi(G_2).
\end{align*}
Since $(p_1^{\alpha_1},\psi(\mathbb{Z}_{p_1^{\alpha_1}}))=(p_1^{\alpha_1}, \frac{p_1^{2\alpha_1+1}+1}{p_1+1})=1$, we deduce that
\begin{align}\label{r2}
\psi(\mathbb{Z}_{p_1^{\alpha_1}})\bigg|\psi\bigg(\xmare{i\in A}{}\mathbb{Z}_{q_i^{\beta_i}} \bigg)-\psi(G_2),
\text{ which further implies that } \psi(\mathbb{Z}_{p_1^{\alpha_1}})\leq\psi\bigg(\xmare{i\in A}{}\mathbb{Z}_{q_i^{\beta_i}} \bigg)-1.
\end{align}
The function $f:[2,\infty)\longrightarrow (0,\infty)$ given by $f(x)=\frac{x^{2\alpha_1+1}+1}{x+1}$ is strictly increasing and $p_1\geq \prod\limits_{i\in A}q_i+1$, so
\begin{align}\label{r3}
\psi(\mathbb{Z}_{p_1^{\alpha_1}})\geq \frac{(\prod\limits_{i\in A}q_i+1)^{2\alpha_1+1}+1}{(\prod\limits_{i\in A}q_i+1)+1}>\prod\limits_{i\in A}q_i^{2\alpha_1}-1.
\end{align}
On the other hand, since $\alpha_1\geq\max\limits_{i\in A}\beta_i$, we have
\begin{align}\label{r4}
\psi\bigg(\xmare{i\in A}{}\mathbb{Z}_{q_i^{\beta_i}} \bigg)-1=\prod\limits_{i\in A}\psi(\mathbb{Z}_{q_i^{\beta_i}})-1<\prod\limits_{i\in A}q_i^{2\beta_i}-1\leq \prod\limits_{i\in A}q_i^{2\alpha_1}-1.
\end{align}
By (\ref{r2}), (\ref{r3}) and (\ref{r4}), we arrive at a contradiction. Therefore, $G$ is not $\psi$-divisible, as desired. 
\hfill\rule{1,5mm}{1,5mm}\\

As a consequence of Proposition 2.4, we highlight the following particular result.\\

\textbf{Corollary 2.5.} \textit{Let $\alpha, \beta$ be integers such that $\alpha\geq 1, \beta\geq 2,$ and let $p$ be an odd prime number. The dihedral groups $D_{2p^{\alpha}}\cong ZM(p^{\alpha}, 2, p^{\alpha}-1)$ and the dicyclic groups $Dic_{4p^{\beta}}\cong ZM(p^{\beta},4,p^{\beta}-1)$ are not $\psi$-divisible.}\\

We note that Proposition 2.4 indicates that $ZM(p,n,r)$ is not $\psi$-divisible, for any odd prime number $p$ and any square-free number $n$. What can be said about the $\psi$-divisibility property if the parameter $m$ of $ZM(m,n,r)$ is a square-free number which is divisible by at least two odd primes? In other words, the following open problem is suggested.\\  

\textbf{Open problem.} Are there any finite non-nilpotent groups of square-free order that are $\psi$-divisible?\\ 

We end this section by mentioning that, using GAP, we saw that all non-nilpotent groups $G$ of square-free order, with $|G|\leq 72000$, are not $\psi$-divisible. In this process, 101021 groups failed to satisfy the $\psi$-divisibility property.

\section{Restricting to normal subgroups}

Since it is not clear if there exist other $\psi$-divisible groups besides the ones indicated by Theorem 1.1, in this section, we choose to work with a less restrictive condition related to the sum of elements order of a finite group. More exactly, we say that a finite group $G$ is \textit{$\psi$-normal divisible} if 
$$\psi(H)|\psi(G), \ \forall \ H\in N(G),$$
where $N(G)$ is the normal subgroup lattice of $G$. The origin of this idea is related to the existence of finite groups having their order equal to the sum of the orders of all their proper subgroups (note that the trivial subgroup is also considered to be proper here). All such groups are isomorphic to $\mathbb{Z}_n$, where $n$ is a perfect number (see Theorem 1 of \cite{7}). However, if one studies finite groups having their order equal to the sum of the orders of all their proper normal subgroups (i.e. the so-called Leinster groups; see \cite{8, 14} for details), more examples can be indicated besides cyclic groups. 

As a consequence of Theorem 1.1, it is clear that all abelian groups of square-free order are $\psi$-normal divisible. \textit{Are there other $\psi$-normal divisible groups?} The answer is clearly yes, since all finite non-abelian simple groups are $\psi$-normal divisible. Moreover, our following result shows that there are infinitely many non-nilpotent, non-simple groups which are $\psi$-normal divisible.\\

\textbf{Theorem 3.1.} \textit{Let $S$ be a finite non-abelian simple group and let $p$ be a prime number such that $(|S|, p)=1$. Then $S\times \mathbb{Z}_p$ is a $\psi$-normal divisible group. In particular, if $p$ is a prime number such that $(60, p)=1$, then $A_5\times \mathbb{Z}_p$ is a $\psi$-normal divisible group.}\\
 
\textbf{Proof.} Let $S$ be a finite non-abelian simple group and let $p$ be a prime number such that $(|S|, p)=1$. Then, the subgroup lattice of $S\times\mathbb{Z}_p$ is decomposable, so the normal subgroups of this group are
$$H\times K, \text{ where } H\in N(S) \text{ and } K\in N(\mathbb{Z}_p).$$
Since both $S$ and $\mathbb{Z}_p$ are $\psi$-normal divisible, we have $\psi(H)\psi(K)|\psi(S)\psi(\mathbb{Z}_p)$. By Lemma 2.1 ii), it follows that
$$\psi(H\times K)|\psi(S\times\mathbb{Z}_p),$$
which completes our proof.
\hfill\rule{1,5mm}{1,5mm}\\

To find other $\psi$-normal divisible groups, the most natural candidates are the finite groups having "few" normal subgroups. In \cite{17} (see Theorem 1.1), the reader can find a classification of finite groups having all normal subgroups of the same order. Such a group $G$ is isomorphic to
\begin{itemize}
\item[--] a simple group;
\item[--] a group having an unique proper normal subgroup;
\item[--] $S\times S$, where $S$ is a simple group;
\item[--] $A_8 \times PSL(3,4)$;
\item[--] $B_n(q)\times C_n(q)$, where $n\geq 3$ and $q$ is an odd prime power.
\end{itemize}
As we previously stated, all finite simple groups are $\psi$-normal divisible. One can use GAP to check that $A_8\times PSL(3,4)$ is not $\psi$-normal divisible. Further, we focus on finite groups having an unique proper normal subgroup, while the other two remaining cases are open. In the same paper, the authors classify the solvable groups having only one proper normal subgroup (see Theorem 1.2) and they indicate the characteristics of non-solvable groups having this property (see Theorem 1.3). In what follows, we are interested in the solvable case.\\

\textbf{Proposition 3.2.} \textit{Let $G$ be a finite solvable group having an unique proper normal subgroup. Then $G$ is not $\psi$-normal divisible.}\\

\textbf{Proof.} Let $G$ be a finite solvable group having only one proper normal subgroup. According to Theorem 1.2 of \cite{17}, we have $G\cong \mathbb{Z}_{p^2}$ or $G\cong \mathbb{Z}_p^{n}\rtimes \mathbb{Z}_q$, where $n\geq 1$ is an integer, $p$ and $q$ are distinct prime numbers, $q|p^n-1$ and $q\nmid p^k-1$, for all $k\in\lbrace 1,2,\ldots, n-1\rbrace$. Since $\mathbb{Z}_{p^2}$ is not a $\psi$-normal divisible group, it remains to study the second case. Assume that $G\cong\mathbb{Z}_{p}^n\rtimes\mathbb{Z}_q$ is $\psi$-normal divisible, where $p,q$ and $n$ have the above mentioned properties. Then,
$$\psi(\mathbb{Z}_{p}^n)|\psi(G)\Longleftrightarrow \psi(\mathbb{Z}_p^n)|\psi(\mathbb{Z}_p^n)+p^n(\psi(\mathbb{Z}_q)-1)\Longleftrightarrow \psi(\mathbb{Z}_p^n)|p^nq(q-1).$$ 
Since $\psi(\mathbb{Z}_p^n)=1+(p^n-1)p$ and $q|p^n-1$, it is easy to show that $(\psi(\mathbb{Z}_p^n),p^n)=(\psi(\mathbb{Z}_p^n),q)=1$ and this leads us to $1+(p^n-1)p|q-1,$ which is false. Hence, $G$ is not $\psi$-normal divisible, as desired.
\hfill\rule{1,5mm}{1,5mm}\\

A study on the non-solvable case is far more difficult to conduct since the characteristics of $G$ depend on the nature of the unique proper normal subgroup $H$. For instance, according to Theorem 1.3 (I) (i) of \cite{17}, if $H$ is solvable and $H\subseteq Z(G)$, then $H\cong\mathbb{Z}_p$, where $p$ is a prime. As an example, for quasisimple groups $G\cong SL(2,q)$, where $q\geq 5$ is an odd prime power, we have $H\cong \mathbb{Z}_2$ and some of these special linear groups are $\psi$-normal divisible, while others are not. For some small values of $q$, using GAP, we remark that $G$ is $\psi$-normal divisible for $q\in \lbrace 5, 7, 9, 17, 19, 27\rbrace$, but the property is not satisfied for $q\in \lbrace 11, 13, 25\rbrace$. 

By Theorem 1.3 (II) (i), if $H$ is non-solvable, then $H$ is isomorphic to a simple non-abelian group $S$ if and only if $G$ is an almost simple group, $S$ is the socle of $G$ and $\frac{G}{S}\cong \mathbb{Z}_p$, where $p$ is a prime number. As an example, for the almost simple groups $G\cong S_n$, where $n\geq 5$, we have $H\cong A_n$. Using GAP, one can check that $S_n$ is not $\psi$-normal divisible for $n\in \lbrace 5,6,\ldots, 12\rbrace$. In this regard, we mention that some useful results on the sum of element orders of symmetric groups were proved in \cite{3} and we suggest the following open problem:\\

\textbf{Open problem.} Let $n\geq 13$ be an integer. Study the $\psi$-normal divisibility property of $S_n$.\\

We end our paper by enumerating the following remarks:
\begin{itemize}
\item[--] Let $G$ be a finite non-nilpotent $\psi$-normal divisible group such that $$|G|\in \lbrace 1,2,\ldots, 1000\rbrace\setminus \lbrace 384, 576, 640, 768, 864, 896, 960\rbrace.$$ Then $G$ is isomorphic to one of the following groups: $A_5, SL(2,5), PSL(2,7), SL(2,7), A_6,$ $A_5\times \mathbb{Z}_7, PSL(2,8), PSL(2,11), A_5\times\mathbb{Z}_{11}, SL(2,9), A_5\times\mathbb{Z}_{13}, SL(2,5)\times \mathbb{Z}_7, PSL(2,7)\times \mathbb{Z}_5$. None of these groups is $\psi$-divisible.
\item[--] There are two non-isomorphic non-nilpotent $\psi$-normal divisible groups having the same order: $PSL(2,11)$ and $A_5\times \mathbb{Z}_{11}$.
\item[--] Let $p$ be a prime number. There are infinitely many $\psi$-normal divisible groups $Q\times \mathbb{Z}_p$, where $Q$ is a quasisimple group such that $(|Q|,p)=1$. In particular $SL(2,5)\times \mathbb{Z}_p$, where $(120,p)=1$, is a $\psi$-normal divisible group. The proof is similar with the one written for Theorem 3.1.
\end{itemize}

\vspace*{3ex}
\small
\hfill
\begin{minipage}[t]{6cm}
Mihai-Silviu Lazorec \\
Faculty of  Mathematics \\
"Al.I. Cuza" University \\
Ia\c si, Romania \\
e-mail: {\tt Mihai.Lazorec@student.uaic.ro}
\end{minipage}
\end{document}